\documentclass[11pt]{article}
\usepackage[margin=1in]{geometry}
\usepackage{amsmath,amssymb,amsthm,mathtools}
\usepackage{enumitem,microtype}
\usepackage{tikz}
\usetikzlibrary{arrows.meta,positioning}
\usepackage[hidelinks,hypertexnames=false]{hyperref}
\newtheorem{theorem}{Theorem}[section]
\newtheorem{lemma}[theorem]{Lemma}
\newtheorem{proposition}[theorem]{Proposition}
\newtheorem{corollary}[theorem]{Corollary}
\newtheorem{assumption}[theorem]{Assumption}
\newtheorem{remark}[theorem]{Remark}
\newcommand{\E}{\mathbb E}\newcommand{\Pp}{\mathbb P}\newcommand{\ns}{\mathrm{ns}}
\title{Iteration-Clock Selection for Euler Schemes of Degenerate Diffusions}
\author{Akihiro Tanaka\thanks{Sumitomo Mitsui Banking Corporation, 1-1-2 Marunouchi, Chiyoda-ku, Tokyo 100-0005, Japan. E-mail: \texttt{tnkaki2000@gmail.com}.}}\date{}
\begin{document}\maketitle

\begin{abstract}
We study Euler--Maruyama approximation of the driftless one-dimensional SDE
$dX_t=\sigma(X_t)\,dW_t$ when $\sigma$ vanishes at the origin and weak uniqueness fails because solutions may remain at zero for an arbitrary time.  The central object is the logarithmic coefficient profile $q(Q)=-\log\sigma(e^{-Q})$.  Under a two-sided logarithmic increment condition, with both indices in $(0,1/2)$, we show that weak-solution selection is governed by a deterministic iteration clock: if $N_q(Q)$ is the number of iterates of $q$ needed to reach a fixed reference depth and $hN_q(\log(1/|x_h|))\to d$, then the Euler scheme started from $x_h$ converges to the fundamental zero-sojourn diffusion delayed by $d$.  This single theorem yields the common-index generalized-power law $d=\theta/|\log\alpha|$ when $h\log\log(1/|x_h|)\to\theta$, an effective translation-number formula for periodically modulated logarithmic profiles, a two-index delay window, and a sharpness construction in which, for one coefficient and one initialization family, the cluster set is the entire delay window.  Thus the macroscopic selection mechanism is the iteration count of the logarithmic coefficient profile; the explicit double-logarithmic law of the homogeneous power case is its common-index closed form.
\end{abstract}



\section{Introduction}

Consider
\[
 dX_t=\sigma(X_t)\,dW_t,\qquad X_0=0,
\]
with a continuous coefficient that is positive away from zero and vanishes at zero. In the regime considered below, the equation admits the Engelbert--Schmidt fundamental zero-sojourn solution together with solutions obtained by inserting a holding period at zero \cite{EngelbertSchmidt1985}. We ask which member of this family is selected by the standard Euler scheme when its deterministic initialization tends to zero as the mesh size tends to zero.

Earlier work \cite[Theorem~9.1]{TaguchiTanaka2019} asserted convergence to the non-sticky diffusion for nonzero initializations approaching a zero of the coefficient without a restriction on their rate of approach. The companion paper \cite[Section~7.6]{TanakaPower2026} gives a counterexample to that unrestricted assertion for the homogeneous power coefficient and establishes the resulting rate-dependent classification into zero-sojourn, delayed, and absorbing limits. Here we identify the structure that persists beyond exact homogeneity: the selector is the iteration count of the logarithmic coefficient profile. For $\sigma(x)=|x|^\alpha$, this gives
\[
 N_q(Q)=\frac{\log Q}{|\log\alpha|}+O(1),
\]
and recovers the power-law selection formula. In the homogeneous setting, the companion paper uses an exactly
rescaled Euler chain, a discrete Tanaka estimate, and a uniform
restart argument. The Tanaka mechanism remains available here,
but the loss of exact scaling requires two-index annular bounds,
weighted estimates for deep returns, and stochastic shadowing of a
coefficient-specific logarithmic iteration.

The problem is related to several strands of the selection and approximation literature. Small-noise selection among nonunique Peano trajectories provides a conceptual analogue in which a perturbation selects one dynamics from a nonunique family \cite{BaficoBaldi1982}. Numerical approximation for SDEs with irregular or non-Lipschitz coefficients has been studied under various uniqueness hypotheses and coefficient regimes; see, among others, \cite{Yamada1978,KanekoNakao1988,Yan2002,BerkaouiBossyDiop2008,GyongyRasonyi2011}. Here the limiting law is itself selected by the numerical scheme from a family of weak solutions. The deterministic holding that arises from this selection mechanism is also distinct from sticky diffusion, which has a different SDE/local-time description \cite{EngelbertPeskir2014}.

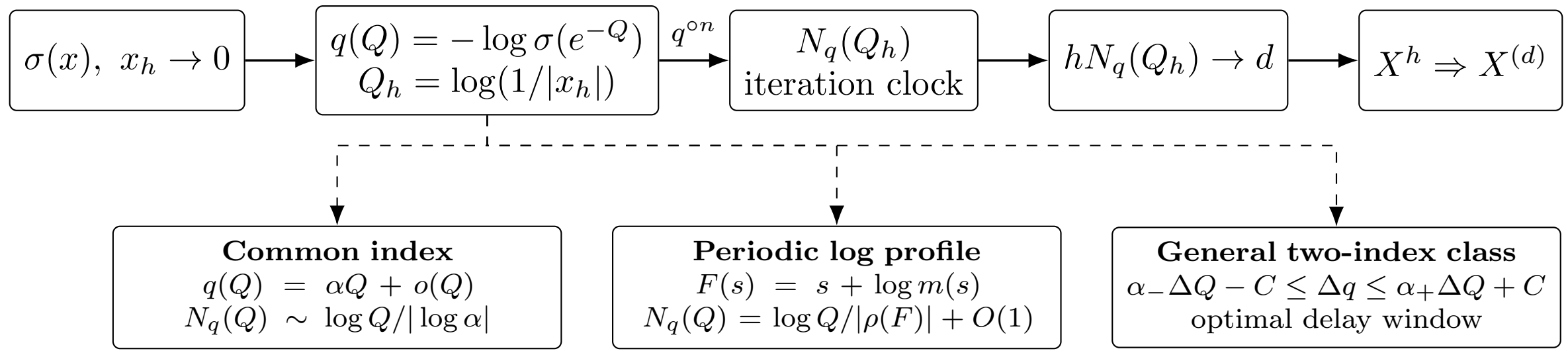
\begin{figure}[t]
\centering
\begin{tikzpicture}[
  font=\small,
  >={Latex[length=2mm]},
  box/.style={draw,rounded corners=2pt,align=center,inner sep=4pt,minimum height=10mm},
  arr/.style={->,line width=.65pt},
  sub/.style={draw,rounded corners=2pt,align=center,inner sep=4pt,text width=42mm,font=\scriptsize}
]
\node[box] (input) {$\sigma(x),\ x_h\to0$};
\node[box, right=7mm of input] (log) {$q(Q)=-\log\sigma(e^{-Q})$\\$Q_h=\log(1/|x_h|)$};
\node[box, right=7mm of log] (clock) {$N_q(Q_h)$\\iteration clock};
\node[box, right=7mm of clock] (delay) {$hN_q(Q_h)\to d$};
\node[box, right=7mm of delay] (law) {$X^h\Rightarrow X^{(d)}$};
\draw[arr] (input)--(log);
\draw[arr] (log)-- node[above,font=\scriptsize] {$q^{\circ n}$} (clock);
\draw[arr] (clock)--(delay);
\draw[arr] (delay)--(law);

\node[sub, below=11mm of log, xshift=-15mm] (power) {\textbf{Common index}\\
$q(Q)=\alpha Q+o(Q)$\\
$N_q(Q)\sim \log Q/|\log\alpha|$};
\node[sub, right=5mm of power] (periodic) {\textbf{Periodic log profile}\\
$F(s)=s+\log m(s)$\\
$N_q(Q)=\log Q/|\rho(F)|+O(1)$};
\node[sub, right=5mm of periodic] (general) {\textbf{General two-index class}\\
$\alpha_-\Delta Q-C\le\Delta q\le\alpha_+\Delta Q+C$\\
optimal delay window};
\draw[dashed,->] (log.south) -- ++(0,-3mm) -| (power.north);
\draw[dashed,->] (log.south) -- ++(0,-3mm) -| (periodic.north);
\draw[dashed,->] (log.south) -- ++(0,-3mm) -| (general.north);
\end{tikzpicture}
\caption{The iteration-clock mechanism.  The coefficient determines a logarithmic dynamical system; its passage-time clock converts microscopic initialization depth into macroscopic holding time.  Constant, periodic, and general two-index profiles correspond respectively to an explicit power-law speed, an effective translation speed, and a delay window.}
\label{fig:clockmechanism}
\end{figure}

The coefficient class is naturally related to $O$-regular variation and Matuszewska indices; see \cite{BGT1987}. We work directly with a two-point increment condition for the logarithmic coefficient profile. Figure~\ref{fig:clockmechanism} summarizes the resulting hierarchy: the power-law constant is the constant-profile case of an effective logarithmic speed; periodic logarithmic profiles admit an effective translation speed; and in the general two-index regime the iteration clock yields a delay window that is sharp in the sense that every point of the window can occur.

\section{Coefficient class and delayed solutions}
Assume throughout that $\sigma$ is continuous, symmetric ($\sigma(-x)=\sigma(x)$) in a neighborhood of zero, positive away from zero, $\sigma(0)=0$, and of at most linear growth. Put
\[
 Q=\log(1/|x|),\qquad q(Q)=-\log\sigma(e^{-Q}).
\]
\begin{assumption}[Two-sided logarithmic increment condition]\label{ass:CI}
For some $Q_\star,C<\infty$ and $0<\alpha_-\le\alpha_+<1/2$,
\begin{equation}
 \alpha_-(Q_2-Q_1)-C\le q(Q_2)-q(Q_1)
 \le\alpha_+(Q_2-Q_1)+C
 \tag{CI}\label{eq:CI}
\end{equation}
for all $Q_2\ge Q_1\ge Q_\star$.
\end{assumption}
No monotonicity of $q$ is assumed.  Condition \eqref{eq:CI} implies
\begin{equation}
 |q(a)-q(b)|\le\alpha_+|a-b|+C,\qquad a,b\ge Q_\star,
 \tag{2.1}\label{eq:coarsecontract}
\end{equation}
and $q(Q)\le\alpha_+Q+C_0$ for large $Q$.

\begin{remark}[Roles of the two thresholds]\label{rem:roles}
The contraction mechanism of the logarithmic iteration uses only that the upper slope bound is strictly below one.  The stronger restriction $\alpha_+<1/2$ ensures local integrability of $\sigma^{-2}$ and the occupation bounds needed to identify the zero-sojourn diffusion after the macroscopic holding period.
\end{remark}

Let $X^{\ns}$ denote the fundamental zero-sojourn solution. For $d<\infty$, let $X^{(d)}$ remain at zero on $[0,d]$ and thereafter follow an independent copy of $X^{\ns}$; put $X^{(\infty)}\equiv0$.

Fix a reference level $Q_{\rm ref}\ge Q_\star$ above the upper affine fixed point $B_+/(1-\alpha_+)$ in Lemma~\ref{lem:detclock}. Such a level will be called admissible. Define the coefficient-specific clock for finite $Q$ by
\begin{equation}
 N_q(Q):=\inf\{n\ge0:q^{\circ n}(Q)\le Q_{\rm ref}\}.
 \tag{2.2}\label{eq:fixedclock}
\end{equation}
We use the extended convention $N_q(+\infty)=+\infty$.

\begin{remark}[Why the clock iterates $q$]\label{rem:qclock}
The Euler depth recursion contains the mesh-dependent additive term $c_h=\tfrac12\log(1/h)$, but the macroscopic number of deep iterations is governed by contraction of the large depth itself. The term $c_h$ determines the intrinsic stopping depth $R_h$ through $R_h-q(R_h)=c_h$; it does not alter the leading iteration count from $Q_h$ down to a mesoscopic level $L_h\gg\log(1/h)$. This is why the coefficient-specific clock is defined from $q^{\circ n}$ rather than from the mesh-dependent map $q+c_h$.
\end{remark}

\begin{lemma}[Deterministic clock bounds]\label{lem:detclock}
There are constants $B_\pm<\infty$ such that, for all sufficiently large $Q$,
\[
 \alpha_-Q-B_-\le q(Q)\le \alpha_+Q+B_+.
 \tag{2.3}\label{eq:affineenvelope}
\]
If $Q_{\rm ref}$ is chosen above the upper affine fixed point
$B_+/(1-\alpha_+)$, then $N_q(Q)<\infty$ for every finite sufficiently
large $Q$.  Moreover, replacing one such fixed reference level by another
changes $N_q(Q)$ by at most a constant independent of $Q$.
Finally, with constants in the $O(1)$ terms allowed to depend on the fixed admissible reference level and on the coefficient constants,
\begin{equation}
 \frac{\log Q}{|\log\alpha_-|}+O(1)
 \le N_q(Q)\le
 \frac{\log Q}{|\log\alpha_+|}+O(1).
 \tag{2.4}\label{eq:clockbounds}
\end{equation}
\end{lemma}
\begin{proof}
The affine envelopes follow from \eqref{eq:CI} with one endpoint fixed.
Iterating the upper envelope gives
\[
 q^{\circ n}(Q)\le
 \alpha_+^nQ+\frac{B_+}{1-\alpha_+},
\]
which proves finite hitting of every admissible reference level. Once
one fixed reference level has been reached, the same affine recursion bounds, uniformly in $Q$,
the number of further iterates needed to reach any
lower admissible fixed level; hence reference-level changes cost only
$O(1)$ iterates. Comparing with the upper and lower affine recursions,
including their fixed-point shifts, gives \eqref{eq:clockbounds}.
\end{proof}

\section{Two-index occupation and zero-sojourn restart}
A fixed rescaled collar means a region $cr_h\le|x|\le Cr_h$ for constants $0<c<C<\infty$ independent of $h$. A sequence of initial points will be called subcritical if it satisfies $h\log(1+\log^+(1/|y_h|))\to0$.

The discrete Tanaka remainder parallels the power-coefficient
argument in \cite{TanakaPower2026}. The occupation and return bounds
below permit coefficients without an exactly rescaled transition
kernel.

For $0<t\le1$ and sufficiently small $x>0$, \eqref{eq:CI} is equivalent to the following bounds with $K_{\rm inc}=e^C$:
\begin{equation}
 K_{\rm inc}^{-1}t^{\alpha_+}\le \frac{\sigma(tx)}{\sigma(x)}
 \le K_{\rm inc}t^{\alpha_-}.
 \tag{3.1}\label{eq:twopower}
\end{equation}
Define
\[
 J(r)=\int_0^r\sigma(a)^{-2}\,da.
\]
\begin{lemma}[Two-index occupation modulus]\label{lem:J}
Under \eqref{eq:twopower},
\[
 \frac{1}{K_{\rm inc}^2(1-2\alpha_-)}\frac r{\sigma(r)^2}
 \le J(r)\le
 \frac{K_{\rm inc}^2}{1-2\alpha_+}\frac r{\sigma(r)^2}.
\]
In particular $J(r)\to0$ as $r\downarrow0$.
\end{lemma}
\begin{proof}
Write $a=rt$ and use \eqref{eq:twopower} in
$J(r)=r\sigma(r)^{-2}\int_0^1[\sigma(rt)/\sigma(r)]^{-2}\,dt$.
\end{proof}

\begin{lemma}[Existence and root-independence of the intrinsic scale]\label{lem:root}
For all sufficiently small \(h\), the equation
\[
 \sqrt h\,\sigma(r)=r
 \tag{3.2}\label{eq:root}
\]
has at least one solution \(r_h\downarrow0\).  If \(r_h^{(1)}\) and
\(r_h^{(2)}\) are two such solutions and
\(R_h^{(i)}=\log(1/r_h^{(i)})\), then
\[
 \bigl|R_h^{(1)}-R_h^{(2)}\bigr|
 \le \frac{C}{1-\alpha_+}.
 \tag{3.3}\label{eq:rootdistance}
\]
Consequently all intrinsic radii are comparable by fixed multiplicative factors, uniformly for small \(h\).
\end{lemma}
\begin{proof}
By the two-sided power comparison, \(\sigma(r)/r\to\infty\) as
\(r\downarrow0\).  At any fixed, sufficiently small radius, the
quantity \(\sqrt h\,\sigma(r)/r\) tends to zero as \(h\downarrow0\).
Continuity therefore gives a root. Since $\sigma(r)\le K(1+r)$ by linear growth, any root satisfies $r\le K\sqrt h/(1-K\sqrt h)$ for small $h$, so all roots tend to zero.

For two roots write \(R_1<R_2\).  Their logarithmic root equations
are
\[
 R_i-q(R_i)=c_h,\qquad c_h=\tfrac12\log(1/h).
\]
Subtracting gives
\[
 R_2-R_1=q(R_2)-q(R_1).
\]
The upper inequality in \eqref{eq:CI} therefore yields
\[
 R_2-R_1\le \alpha_+(R_2-R_1)+C,
\]
which is exactly \eqref{eq:rootdistance}.  Exponentiation gives the
multiplicative comparison of the corresponding radii.
\end{proof}

\begin{lemma}[Intrinsic-scale occupation estimates]\label{lem:intrinsic}
Let \(r_h\downarrow0\) satisfy \(\sqrt h\,\sigma(r_h)=r_h\), and set
\[
a_h(u)=\frac{\sigma(r_h|u|)}{\sigma(r_h)}.
\]
Then
\[
J(r_h)\asymp \frac{h}{r_h}.
\tag{3.4}\label{eq:Jrh}
\]
Moreover, on every fixed collar \(0<R<M<\infty\),
\[
0<c_{R,M}\le a_h(u)\le C_{R,M}<\infty,
\qquad R\le |u|\le M,
\tag{3.5}\label{eq:collar}
\]
uniformly for small \(h\).
\end{lemma}
\begin{proof}
The first assertion follows from Lemma~\ref{lem:J} and
\(\sigma(r_h)^2=r_h^2/h\).  The second follows directly from
\eqref{eq:twopower}, applied separately when \(|u|\le1\) and
\(|u|\ge1\).
\end{proof}

\begin{lemma}[Uniform logarithmic shocks]\label{lem:logshock}
For every fixed $\delta_0<\infty$ there are constants $C,c>0$ such
that, for $Z\sim N(0,1)$, $|\delta|\le\delta_0$, and
$Y_\delta=-\log|Z+\delta|$,
\[
 \Pp(Y_\delta>y)\le Ce^{-y},\qquad
 \Pp(Y_\delta<-y)\le C\exp\{-ce^{2y}\},\qquad y\ge0.
 \tag{3.6}\label{eq:uniformshock}
\]
In particular $\sup_{|\delta|\le\delta_0}\E(Y_\delta^+)<\infty$.
\end{lemma}
\begin{proof}
The first event is $|Z+\delta|<e^{-y}$ and is bounded by the length
$2e^{-y}$ times the supremum of the Gaussian density.  For the second,
$|Z+\delta|>e^y$ implies $|Z|>e^y-\delta_0$; the Gaussian tail gives
the stated bound after enlarging the constants on a bounded range of
$y$.
\end{proof}

\begin{lemma}[Foster--Lyapunov descent to the intrinsic scale]\label{lem:lyapdescent}
Let $r_h$ be an intrinsic radius and $R_h=\log(1/r_h)$.  There exist fixed $B,C<\infty$ such that, whenever
$Q_k\ge R_h+B$, the logarithmic depth satisfies
\[
 Q_{k+1}=q(Q_k)+c_h+Y_{k+1}^h,
 \qquad c_h=\tfrac12\log(1/h),
\]
where $Y_{k+1}^h=-\log|Z_{k+1}+\delta_{h,k}|$ is the logarithmic shock arising from the Euler step.  With
\[
 \tau_B=\inf\{n\ge0:Q_n\le R_h+B\},
\]
\[
 \E_Q\tau_B\le C\{1+\log(1+(Q-R_h)^+)\}.
 \tag{3.7}\label{eq:lyapdescent}
\]
The constants are uniform for small $h$. The exit position
$X^h_{h\tau_B}/r_h$ is also tight uniformly for initial depths
$Q\ge R_h+B$.
\end{lemma}
\begin{proof}
By the intrinsic identity $q(R_h)+c_h=R_h$ and the upper increment
bound, for $U=Q-R_h\ge0$,
\[
 q(R_h+U)+c_h-R_h\le\alpha_+U+C_0.
\]
Moreover the Euler shift satisfies
\[
 |\delta_{h,k}|=\exp\{-[Q-q(Q)-c_h]\}
 \le \exp\{-(1-\alpha_+)U+C_0\},
\]
so it is uniformly bounded on $U\ge B$.  Lemma~\ref{lem:logshock} therefore applies uniformly.  Choose $\lambda\in(\alpha_+,1)$ and enlarge $B$ so that the right-hand
side is at most $\lambda U$ for $U\ge B$.  Hence, before $\tau_B$,
\[
 U_{k+1}^+\le \lambda U_k+(Y_{k+1}^h)^+.
\]
Lemma~\ref{lem:logshock} gives
$\sup\E[(Y_{k+1}^h)^+\mid\mathcal F_k]\le M$.  With
$V(u)=\log(1+u)$, Jensen's inequality yields
\[
 \E[V(U_{k+1}^+)\mid\mathcal F_k]
 \le\log(1+\lambda U_k+M).
\]
After one further enlargement of $B$ there is $c_0>0$ such that the
last expression is at most $V(U_k)-c_0$ for $U_k\ge B$.  Optional
stopping at $\tau_B\wedge n$ gives
$c_0\E(\tau_B\wedge n)\le V(U_0)+C$; monotone convergence proves
\eqref{eq:lyapdescent}.

For the exit position put $\kappa=e^{-B}$.  As long as $k<\tau_B$,
$|X^h_{kh}|<\kappa r_h$ and the two-index bounds give
$\sqrt h\,\sigma(X^h_{kh})\le Cr_h$.  If $Y=x+sZ$ with
$|x|\le\kappa r_h$, $0<s\le Cr_h$ and $Z\sim N(0,1)$, the Gaussian
survival function satisfies
\[
 \frac{\overline\Phi(a+t)}{\overline\Phi(a)}
 \le e^{-t^2/2},\qquad a,t\ge0.
\]
Indeed, after translating the integral defining $\overline\Phi(a+t)$
by $t$, its integrand is at most $e^{-t^2/2}$ times the integrand
defining $\overline\Phi(a)$.
For $\mu=x/r_h\in[-\kappa,\kappa]$ and $v=s/r_h\in(0,C]$,
apply the ratio bound separately to the positive and negative Gaussian
tails, with thresholds $(\kappa\mp\mu)/v\ge0$ and increment
$(K-\kappa)/v$.  Thus, for every $K>\kappa$,
\[
 \Pp(|Y|>Kr_h)\le\eta(K)\Pp(|Y|>\kappa r_h),\qquad
 \eta(K)=\exp\left\{-\frac{(K-\kappa)^2}{2C^2}\right\}.
\]
Summing over the disjoint first-exit events and
conditioning at the preceding grid time gives
\[
 \Pp(|X^h_{h\tau_B}|>Kr_h)
 \le\eta(K)\Pp(\tau_B<\infty)\le\eta(K).
\]
This yields uniform tightness of the rescaled exit position.
\end{proof}

\begin{lemma}[Deep-side subcritical entrance]\label{lem:subcriticalescape}
Let $y_h\to0$, $y_h\ne0$, satisfy
\[
 hH(y_h)\to0,\qquad H(y)=\log\!\left(1+\log^+\frac1{|y|}\right),
 \tag{3.8}\label{eq:subcrit}
\]
and suppose in addition that its depth $Q_h=\log(1/|y_h|)$ satisfies
$Q_h\ge R_h+B$ for a fixed sufficiently large $B$.  If
\[
 \tau_h^B=\inf\{kh:-\log|X^h_{kh}|\le R_h+B\},
\]
then $\tau_h^B\to0$ in probability and
$\E(\tau_h^B\wedge T)\to0$ for every fixed $T$.  The rescaled exit position is tight and lies outside the fixed inner radius $e^{-B}$.
\end{lemma}
\begin{proof}
Choose $L_h\uparrow\infty$ with $\log(1/h)=o(L_h)$ and
$h\log L_h\to0$. If $Q_h\le L_h$, the shadowing stage below is empty and we proceed directly to the Foster--Lyapunov descent of Lemma~\ref{lem:lyapdescent}. Otherwise, above $L_h$, compare the stochastic depth with
$D_{k+1}=q(D_k)$ only until the first crossing of $L_h$ by either
orbit.  Write $\nu_h=\inf\{k\ge0:Q_k\le L_h\}$. Lemma~\ref{lem:logshock}
and a union bound give, on every fixed physical horizon and for any
fixed $A>2$,
\[
 \Pp\!\left(\max_{1\le k\le\lfloor K/h\rfloor\wedge\nu_h}
 |Y_k^h|>A\log(1/h)\right)\longrightarrow0,
\]
where an empty maximum is zero. Indeed the positive-shock tail is bounded by
$C(K/h)e^{-A\log(1/h)}$, and the negative-shock tail is
smaller. Let $D_{k+1}=q(D_k)$, $D_0=Q_h$,
and stop both depth processes at their first crossing of $L_h$.
The two-point increment bound for $q$, together with the logarithmic
recursion, yields on this stopped interval
\[
 |Q_{k+1}-D_{k+1}|
 \le \alpha_+|Q_k-D_k|+C+\tfrac12\log(1/h)+|Y_{k+1}^h|.
\]
Iteration and the preceding maximal-shock bound give
\[
 \max_{k\le K/h,\;k<\tau_{L_h}}|Q_k-D_k|
 =O_{\Pp}(\log(1/h))=o_{\Pp}(L_h),
\]
because $\log(1/h)=o(L_h)$.  Lemma~\ref{lem:detclock} and $hH(y_h)\to0$ show that
the deterministic crossing time of $L_h$ is $o(1)$ in physical time;
the same is true for the Euler crossing time.  Indeed the two
stopped depth processes differ by at most a deterministic
$b_h=C_A\log(1/h)=o(L_h)$ on the preceding event, and the shocks
are bounded by $s_h=A\log(1/h)$. For $h$ small,
\[
 \alpha_+(L_h+b_h)+C_0+\tfrac12\log(1/h)+s_h<L_h-b_h.
\]
Once either depth reaches $L_h$, the other therefore crosses within
one more iteration. The stochastic crossing depth belongs to
$[cL_h,CL_h]$ with probability tending to one by the two-sided
affine bounds.

From that point to $R_h+B$, Lemma~\ref{lem:lyapdescent} gives an
expected additional number of steps $O(1+\log L_h)$, hence physical
time $o(1)$.  Before this entrance the grid positions are at most
$e^{-B}r_h$ in absolute value and their one-step diffusion amplitudes
are at most $Cr_h$. The maximal normalized Brownian increment over
$T/h$ cells is $O_{\Pp}(\sqrt{\log(1/h)})$. Thus the continuous
path before entrance converges uniformly to zero in probability,
and Lemma~\ref{lem:lyapdescent} provides tightness of its
rescaled exit position.
\end{proof}

\begin{lemma}[Discrete Tanaka boundary estimate]\label{lem:tanaka}
Fix a localization radius $M<\infty$ and a horizon $T<\infty$.  There are constants
$h_0,c,C>0$ such that the following holds.  Let $a>0$ be sufficiently
small, let $0<h\le h_0$, and assume, for one fixed $\varepsilon_0\in(0,1)$, that
\[
 c\sqrt h\,\sigma(a)\le \varepsilon_0 a.
\]
Assume also that the boundary layer
\[
 \mathcal B_{a,h}=\left\{x:\bigl||x|-a\bigr|
 \le c\sqrt h\,\sigma(a)\right\}
\]
is contained in $(-M,M)$.  Then, for the Euler chain stopped on
leaving $(-M,M)$,
\[
 h\E\sum_{k<T/h}\mathbf1_{\{X^h_{kh}\in\mathcal B_{a,h}\}}
 \le C_{T,M,\varepsilon_0}\frac{\sqrt h}{\sigma(a)}.
 \tag{3.9}\label{eq:tanakaocc}
\]
The constants are uniform in $a$ and $h$ under the displayed condition (with fixed $\varepsilon_0$). The layer lies in a fixed multiplicative neighborhood of $a$.
\end{lemma}
\begin{proof}
Set $f_a(x)=(|x|-a)^+$ and let $f'_{a,-}$ be its left derivative.
For $\Delta_k=X^h_{(k+1)h}-X^h_{kh}$ define the convexity remainder
\[
 R_a(x,z)=f_a(x+z)-f_a(x)-f'_{a,-}(x)z\ge0.
\]
At deterministic grid times,
\[
 f_a(X^h_{Nh})-f_a(X^h_0)
 =\sum_{k<N}f'_{a,-}(X^h_{kh})\Delta_k
  +\sum_{k<N}R_a(X^h_{kh},\Delta_k).
 \tag{3.10}\label{eq:disctanaka}
\]
After stopping at the localization boundary, the first sum is a true
martingale.  If $\tau$ is the stopped terminal grid time, then
\[
 \E\sum_{k<\tau}R_a(X^h_{kh},\Delta_k)
 =\E f_a(X^h_\tau)-f_a(X^h_0)
 \le \E f_a(X^h_\tau).
\]
For $k<\tau$, the coefficient $\sigma(X^h_{kh})$ is bounded by
$K_M$.  Hence
\[
 \E f_a(X^h_{\tau h})
 \le M+K_M\E\max_{k<T/h}|W_{(k+1)h}-W_{kh}|
 \le M+C_{M,T}\sqrt{h\log(2+T/h)}.
\]
In particular the remainder sum has an upper bound independent of
$a$ and $h$.

Condition on $X^h_{kh}=x$ and put $s=\sqrt h\,\sigma(x)$.  Near the
positive kink $a$, the contribution of $g_a(x)=(x-a)^+$ is explicit:
if $z=|x-a|/s$ and $\overline\Phi(z)=1-\Phi(z)$, then
\[
 \E R_a(x,sZ)\ge
 s\{\varphi(z)-z\,\overline\Phi(z)\},
 \tag{3.11}\label{eq:gaussremainder}
\]
up to the nonnegative contribution of the opposite kink.  This is the identity
$\E(Z-z)^+=\varphi(z)-z\overline\Phi(z)$; the case $x>a$, involving the corresponding $(a-x-sZ)^+$ term, gives the same expression.  The same formula holds near $-a$ by symmetry.  On $\mathcal B_{a,h}$ the
two-index comparison gives $\sigma(x)\asymp\sigma(a)$; after choosing
the layer constant $c$ fixed, $z$ is uniformly bounded.  Since
$\varphi(z)-z\overline\Phi(z)$ is strictly positive on bounded
intervals,
\[
 \E[R_a(x,\sqrt h\sigma(x)Z)\mid X^h_{kh}=x]
 \ge c_0\sqrt h\,\sigma(a),
 \qquad x\in\mathcal B_{a,h}.
 \tag{3.12}\label{eq:remainderlower}
\]
Summing \eqref{eq:remainderlower} in \eqref{eq:disctanaka} gives
\eqref{eq:tanakaocc}.  
\end{proof}

\begin{lemma}[Mesoscopic deterministic-grid occupation]\label{lem:gridocc}
Fix $R>1$, $T<\infty$, and a localization radius $M$.  For every
sufficiently small $\varepsilon<M$,
\[
 \limsup_{h\downarrow0}
 \E\left[h\sum_{k<T/h}
 \mathbf1_{\{Rr_h\le |X^h_{kh}|<\varepsilon\}}\right]
 \le C_{T,M,R}\int_0^{C\varepsilon}\frac{da}{\sigma(a)^2}.
 \tag{3.13}\label{eq:annocc}
\]
In particular the right-hand side tends to zero as
$\varepsilon\downarrow0$.
\end{lemma}
\begin{proof}
On a dyadic annulus
$[a,2a]$, the two-index bounds imply $\sigma(x)\asymp\sigma(a)$.
For the annuli used here, $a\ge Rr_h$.  The intrinsic identity and the two-index comparison imply
$\sqrt h\,\sigma(a)/a\le C_R$, with $C_R\downarrow0$ as $R\uparrow\infty$; after fixing the outer collar parameter $R$ (and, if necessary, decreasing the fixed Tanaka-layer constant), the admissibility condition in Lemma~\ref{lem:tanaka} holds uniformly in $h$ on all these annuli.
Starting from $a_0=a$, choose successively
\[
 a_{m+1}=a_m+\lambda_0\sqrt h\,\sigma(a_m)
\]
with a fixed small $\lambda_0>0$, and stop after crossing $2a$.
By the two-index comparison, $\sigma(a_m)$ and $\sigma(a_{m+1})$
are comparable uniformly on the annulus. Choosing $\lambda_0$
smaller than the fixed Tanaka layer-width constant, after allowing
for this comparison, makes the layers centered at successive $a_m$
cover $[a,2a]$. Their overlap multiplicity is bounded by a constant
independent of $h$: throughout the annulus all layer widths are
comparable to $\sqrt h\sigma(a)$, and the spacing between centers
is comparable to $\lambda_0\sqrt h\sigma(a)$.
Lemma~\ref{lem:tanaka} therefore gives
\[
 h\E\sum_{k<T/h}\mathbf1_{\{a\le|X^h_{kh}|<2a\}}
 \le C_{T,M,R}\sum_m\frac{\sqrt h}{\sigma(a_m)}.
 \tag{3.14}\label{eq:layersum}
\]
Because $a_{m+1}-a_m=\lambda_0\sqrt h\sigma(a_m)$, the last sum is a
Riemann sum for the speed density:
\[
 \sum_m\frac{\sqrt h}{\sigma(a_m)}
 =\lambda_0^{-1}\sum_m
 \frac{a_{m+1}-a_m}{\sigma(a_m)^2}
 \le C\int_{a/2}^{4a}\frac{du}{\sigma(u)^2}.
 \tag{3.15}\label{eq:speedriemann}
\]
Only the two endpoint intervals extend beyond the nominal annulus;
their contributions are included in the enlarged integration interval
$[a/2,4a]$.
Summing over dyadic annuli from $Rr_h$ to $\varepsilon$ gives
\eqref{eq:annocc}; bounded overlap changes only the constant.  The
integral tends to zero by Lemma~\ref{lem:J}.
\end{proof}

\begin{lemma}[Occupation near the degenerate point]\label{lem:deepreturn}
Fix a localization radius $M$ and $T<\infty$.  Let $B$ be the
constant in Lemma~\ref{lem:lyapdescent}, and put $\kappa=e^{-B}$.
For any sequence of initial points $y_h\to0$ with
$|y_h|\ge\kappa r_h$, and uniformly for initial points in a fixed
rescaled collar bounded away from zero,
\begin{equation}
 \lim_{\varepsilon\downarrow0}\limsup_{h\downarrow0}
 \E\int_0^T\mathbf1_{\{|X_t^h|<\varepsilon\}}\,dt=0.
 \tag{3.16}\label{eq:contocc}
\end{equation}
The same assertion holds for the chain stopped upon leaving $(-M,M)$,
with constants depending on $M$ and $T$.
\end{lemma}
\begin{proof}
We first work at grid times and stop the chain at its first grid exit
from $(-M,M)$.  Choose $R$ larger than $2$ and the fixed outer-collar threshold required in Lemma~\ref{lem:gridocc} and in the Gaussian far-return estimate, and
$0<\varepsilon<\varepsilon_*$, where $\varepsilon_*$ is inside the
neighborhood on which \eqref{eq:coarsecontract} holds.  For small $h$,
the grid occupation of $\{|x|<\varepsilon\}$ is divided into
\[
 \{|x|<\kappa r_h\},\qquad
 \{\kappa r_h\le|x|<Rr_h\},\qquad
 \{Rr_h\le|x|<\varepsilon\}.
\]
The last region satisfies Lemma~\ref{lem:gridocc}.  On the fixed
rescaled annulus $[\kappa r_h,Rr_h]$, the two-index comparison and
$\sqrt h\,\sigma(r_h)=r_h$ imply
\[
 0<c_{\kappa,R}\le
 \frac{\sigma(r_h|u|)}{\sigma(r_h)}
 \le C_{\kappa,R},\qquad \kappa\le|u|\le R.
\]
Choose the layer-width constant in Lemma~\ref{lem:tanaka} sufficiently
small in terms of $\kappa$ and $R$.  The fixed rescaled annulus can
then be covered by a number of admissible Tanaka layers bounded
independently of $h$.  Each layer has width comparable to $r_h$ and
satisfies $h\E N_{\rm layer}\le C_{\kappa,R,T}h/r_h$.  Thus
\begin{equation}
 \E\left[h\sum_{k<T/h}
 \mathbf1_{\{\kappa r_h\le|X^h_{kh}|<Rr_h\}}\right]
 \le C_{\kappa,R,T}\frac h{r_h}.
 \tag{3.17}\label{eq:collarocc}
\end{equation}

It remains to account for grid times with $|X^h_{kh}|<\kappa r_h$.
We group consecutive such times into disjoint excursions, each
starting at the first grid entry from the complement and ending on
its first subsequent grid exit.  No initial excursion is present
under the stated initialization condition.  If $Y$ is the landing
point of an entry, put
$U(Y)=\log(\kappa r_h/|Y|)>0$.
Lemma~\ref{lem:lyapdescent}, applied to
$Q(Y)=R_h+B+U(Y)$, gives the conditional mean physical duration of
that excursion bounded by
\begin{equation}
 \E_Y\tau_{\mathrm{exit}}
 \le Ch\{1+\log(1+U(Y))\}.
 \tag{3.18}\label{eq:excursioncost}
\end{equation}
The bound also holds for the duration truncated at $T$.

We distinguish entries according to the preceding grid position
$x=X^h_{kh}$.  For $\kappa r_h\le|x|\le Rr_h$, let
$Y=x+\sqrt h\sigma(x)Z$, $Z\sim N(0,1)$.  On this fixed rescaled
annulus, the Gaussian transition density is bounded by $C/r_h$.
Consequently, for all $s\ge0$,
\[
 \Pp_x\{|Y|<\kappa r_h e^{-s}\}\le Ce^{-s}.
\]
Since $g(s)=1+\log(1+s)$ has $g(0)=1$ and
$g'(s)=(1+s)^{-1}$, integrating the preceding tail estimate yields
\[
 \E_x[\mathbf1_{\{|Y|<\kappa r_h\}}g(U(Y))]
 =\Pp_x(|Y|<\kappa r_h)
 +\int_0^\infty\frac{\Pp_x(|Y|<\kappa r_h e^{-s})}{1+s}\,ds.
\]
In particular,
\begin{equation}
 \sup_{\kappa r_h\le|x|\le Rr_h}
 \E_x\left[\mathbf1_{\{|Y|<\kappa r_h\}}
 \{1+\log(1+U(Y))\}\right]\le C_{\kappa,R}.
 \tag{3.19}\label{eq:nearlanding}
\end{equation}
Let $\mathcal I_{\mathrm{near}}$ denote the grid indices $k<T/h$
for which the stopped chain is in $[\kappa r_h,Rr_h]$ in absolute
value and the next grid state enters $(-\kappa r_h,\kappa r_h)$.
For each such index, the expected remaining physical time in that
deep excursion, conditional on the landing point, is bounded by
\eqref{eq:excursioncost}. Excursions are disjoint, and truncating
one at $T$ can only decrease its contribution. Conditioning first
on the current grid state and then on the landing point gives
\[
\begin{aligned}
&\E\sum_{k\in\mathcal I_{\mathrm{near}}}
  (\text{physical duration of the resulting deep excursion})\\
&\qquad\le Ch\sum_{k<T/h}
  \E\!\left[\mathbf1_{\{\kappa r_h\le|X^h_{kh}|\le Rr_h\}}
  \E_{X^h_{kh}}\!\left[
  \mathbf1_{\{|Y|<\kappa r_h\}}
   \{1+\log(1+U(Y))\}\right]\right]\\
&\qquad\le C_{\kappa,R,T}h/r_h,
\end{aligned}
\]
where the last inequality uses \eqref{eq:nearlanding}
and \eqref{eq:collarocc}. Thus the deep occupation generated
by these near-annulus entries vanishes.

For predecessor states farther away, write $x=r_hu$ with
$R\le|u|\le\varepsilon_*/r_h$.  The two-index comparison implies
\[
 \frac{|x|}{\sqrt h\sigma(x)}\ge
 c|u|^{1-\alpha_+},\qquad
 \frac{r_h}{\sqrt h\sigma(x)}\le C|u|^{-\alpha_-}.
\]
A Gaussian density estimate on the interval
$(-\kappa r_h e^{-s},\kappa r_h e^{-s})$ therefore yields, for
$s\ge0$ and a constant $c>0$,
\[
 \Pp_x\{|Y|<\kappa r_h e^{-s}\}
 \le C\exp\{-c|u|^{2(1-\alpha_+)}\}\,e^{-s}.
\]
Integrating the logarithmic weight by the same tail-integration
identity, and using $\int_0^\infty e^{-s}(1+s)^{-1}\,ds<\infty$,
we obtain
\begin{equation}
 \E_x\left[\mathbf1_{\{|Y|<\kappa r_h\}}
 \{1+\log(1+U(Y))\}\right]
 \le C\exp\{-c|u|^{2(1-\alpha_+)}\}.
 \tag{3.20}\label{eq:farlanding}
\end{equation}
Set $a_j=2^jRr_h$ and decompose the local far region into dyadic
annuli $[a_j,2a_j]$.  Equations
\eqref{eq:layersum}--\eqref{eq:speedriemann} show that
\[
 \E N_{j,h}(T)\le \frac{C_T}{h}
 \int_{a_j/2}^{4a_j}\frac{da}{\sigma(a)^2}.
\]
The same conditioning-and-summation argument, with
\eqref{eq:farlanding} in place of \eqref{eq:nearlanding},
bounds the total expected duration of deep excursions entered
from these annuli by
\[
 C_T\sum_j e^{-c(2^jR)^{2(1-\alpha_+)}}
 \int_{a_j/2}^{4a_j}\frac{da}{\sigma(a)^2}.
\]
For annuli contained in the local coefficient neighborhood,
\eqref{eq:coarsecontract} and the intrinsic identity give
\[
 \int_{a_j/2}^{4a_j}\frac{da}{\sigma(a)^2}
 \le C\frac h{r_h}(2^jR)^{1-2\alpha_-}.
\]
Hence their contribution is bounded by
\begin{equation}
 C_T\frac h{r_h}\sum_{j\ge0}
 (2^jR)^{1-2\alpha_-}
 e^{-c(2^jR)^{2(1-\alpha_+)}}
 =O_T(h/r_h).
 \tag{3.21}\label{eq:farcost}
\end{equation}
If the predecessor lies outside the local coefficient neighborhood
but inside $(-M,M)$, choose the neighborhood endpoint
$\varepsilon_*>0$ as above.  On the compact set
$\{x:\varepsilon_*\le |x|\le M\}$, the coefficient is bounded
above and away from zero.  The shifted-Gaussian density at zero is
therefore bounded by $C_M h^{-1/2}e^{-c_M/h}$, uniformly over these
predecessor states.  Integrating over the landing interval of width
$2\kappa r_h e^{-s}$, and then integrating the logarithmic weight,
shows that the contribution per grid transition is at most
$C_M r_h h^{-1/2}e^{-c_M/h}$.  At most $T/h+1$ grid transitions are
possible, so the total expected deep duration from this region is
bounded by $C_{M,T}r_h h^{-1/2}e^{-c_M/h}$ and tends to zero.  Since $h/r_h\asymp J(r_h)\to0$, the expected grid
occupation of the deep core tends to zero.

Combining the three spatial regions gives, for fixed $\varepsilon$,
\[
 \limsup_{h\downarrow0}\E\left[h\sum_{k<T/h}
 \mathbf1_{\{|X^h_{kh}|<\varepsilon\}}\right]
 \le C_{T,M,R}\int_0^{C\varepsilon}\frac{da}{\sigma(a)^2}.
\]
The right-hand side tends to zero with $\varepsilon$.

For the continuous Euler interpolation, let
$\Delta_h(T)=\max_{k<T/h}\sup_{0\le s\le h}
 |X^h_{kh+s}-X^h_{kh}|$.  Under localization,
$\sigma$ is bounded, and the Brownian maximal inequality and a
union bound give $\Delta_h(T)\to0$ in probability.
Consequently, for each fixed $\varepsilon>0$,
\[
 \E\int_0^T\mathbf1_{\{|X^h_t|<\varepsilon\}}\,dt
 \le h\E\sum_{k<T/h}
 \mathbf1_{\{|X^h_{kh}|<2\varepsilon\}}
 +T\Pp(\Delta_h(T)>\varepsilon)+O(h).
\]
Taking $h\downarrow0$ and then $\varepsilon\downarrow0$
proves \eqref{eq:contocc} for the localized chain.  Finally,
linear growth and the standard BDG--Gronwall estimate give
$\sup_{h\le h_0}\Pp(\sup_{t\le T}|X_t^h|\ge M)\to0$
as $M\to\infty$ (for initial points tending to zero).  Since the
occupation integral is bounded by $T$, removing localization gives
the stated assertion.
\end{proof}

\begin{proposition}[Subcritical zero-sojourn identification]\label{prop:restart}
Let $y_h\to0$, $y_h\ne0$, satisfy $hH(y_h)\to0$.  Then the Euler
interpolation started from $y_h$ converges to $X^{\ns}$.  The
convergence is uniform when the restart points range over a fixed
rescaled collar.
\end{proposition}
\begin{proof}
There are two cases.  Along a deep-side
subsequence, $|y_h|<e^{-B}r_h$, Lemma~\ref{lem:subcriticalescape}
brings the chain to a fixed intrinsic collar in $o_{\Pp}(1)$ physical
time while the preceding physical path vanishes.  Along an outer-side
subsequence, $|y_h|\ge e^{-B}r_h$, no collar entrance is imposed:
Lemma~\ref{lem:deepreturn} applies from time zero and gives vanishing
occupation near the origin directly.  Thus every sequence admits a
subsequence on which one of these two arguments applies.

After localization, linear growth, BDG, and Gronwall give tightness of
the Euler interpolations in $C([0,T])$, and the grid interpolation
error tends to zero in probability.  Every subsequential limit solves
the martingale problem for
\[
 \mathcal Lf(x)=\tfrac12\sigma(x)^2f''(x).
 \tag{3.22}\label{eq:generator}
\]
The occupation bounds imply
\[
 \int_0^T\mathbf1_{\{X_t=0\}}\,dt=0\quad\text{a.s.}
 \tag{3.23}\label{eq:zerosojourn}
\]
To justify the passage to the limit, for fixed $\varepsilon>0$
the functional $\omega\mapsto\int_0^T
\mathbf1_{\{|\omega(t)|<\varepsilon\}}dt$ is lower semicontinuous
on $C([0,T])$.  Portmanteau and the bounds from
Lemma~\ref{lem:deepreturn}, followed by $\varepsilon\downarrow0$,
therefore give \eqref{eq:zerosojourn}. For deep-side initializations,
the initial entrance interval contributes $o(1)$ in expected
truncated duration by Lemma~\ref{lem:subcriticalescape}, and its
continuous path is negligible.
Since $\sigma$ is continuous, every continuous limit of the stopped
Euler martingales solves the corresponding local martingale problem,
and hence is a continuous local martingale with quadratic variation
$\int_0^t\sigma(X_s)^2\,ds$.  The standard Brownian representation
of a continuous local martingale (enlarging the probability space if
needed on the set where $\sigma=0$) gives a weak solution of the
SDE. By Lemma~\ref{lem:J}, $\sigma^{-2}$ is locally integrable.
The zero-sojourn law is uniquely determined by the natural-scale time change: if $B$ is Brownian motion started at zero, the inverse of $A_u=\int_0^u\sigma(B_s)^{-2}\,ds$ (with arbitrary values at $B_s=0$) determines the fundamental diffusion, and zero sojourn excludes any additional holding time at the origin. The local integrability from Lemma~\ref{lem:J} makes this time change finite on compact Brownian-time intervals. This is the fundamental-solution construction of Engelbert and Schmidt \cite{EngelbertSchmidt1985}, and it identifies the limit as $X^{\ns}$.

Uniformity on a fixed rescaled collar follows by sequential
contradiction.  Any violating sequence of collar points has a
subsequence whose rescaled locations converge in the compact collar;
the same tightness and occupation estimates then force convergence to
$X^{\ns}$.
\end{proof}

\begin{lemma}[Residual descent after clock completion]\label{lem:residual}
Let $L_h\uparrow\infty$ satisfy $\log(1/h)=o(L_h)$ and
$h\log L_h\to0$.  If the Euler depth starts in
$[c_0L_h,C_0L_h]$, the additional physical time needed to reach a
fixed enlarged intrinsic collar is $o_{\Pp}(1)$; its truncated
expectation is $o(1)$ on every fixed horizon.  During this residual
phase the physical path converges uniformly to zero in probability,
and the rescaled exit state is tight.
\end{lemma}
\begin{proof}
The intrinsic identity is $q(R_h)+c_h=R_h$ and
$R_h=O(\log(1/h))$.  The initial excess above $R_h$ is $O(L_h)$.
Lemma~\ref{lem:lyapdescent} therefore bounds the expected number of
additional steps by $C\{1+\log L_h\}$, whose physical duration is
$o(1)$.  Before that exit the grid positions are bounded by
$e^{-B}r_h$ and their diffusion amplitudes by $Cr_h$.
Brownian maximal estimates therefore imply uniform disappearance
of the continuous path before exit in probability, while
Lemma~\ref{lem:lyapdescent} gives tightness of the rescaled entrance
position.
\end{proof}

\section{The iteration-clock theorem}
Let $x_h\to0$ be deterministic and
\[
 X_t^h=x_h+\int_0^t\sigma(X^h_{h\lfloor s/h\rfloor})\,dW_s.
\]
For $x_h\ne0$ put $Q_h=\log(1/|x_h|)$.  If $x_h=0$, use the extended convention $Q_h=+\infty$ and $N_q(Q_h)=+\infty$; in that case the Euler scheme is identically zero because $\sigma(0)=0$.
\begin{theorem}[Iteration-clock selection]\label{thm:master}
Under Assumption~\ref{ass:CI}, if
\[
 hN_q(Q_h)\to d\in[0,\infty),
\]
then $X^h\Rightarrow X^{(d)}$ in $C_{\rm loc}([0,\infty))$.  If $hN_q(Q_h)\to\infty$, then $X^h\to0$ in probability in $C_{\rm loc}([0,\infty))$.
\end{theorem}
\begin{proof}
On any subsequence with $x_h=0$, the Euler path is identically zero.
Such a subsequence is compatible only with a divergent clock; in that
case the asserted limit is immediate. We may therefore work on
subsequences with $x_h\ne0$.

Choose an auxiliary cutoff $L_h\uparrow\infty$ with
\begin{equation}
 \log(1/h)=o(L_h),\qquad h\log L_h\to0.
 \tag{4.1}\label{eq:Lh}
\end{equation}
The number of $q$-iterations between $L_h$ and the fixed level $Q_{\rm ref}$ is $O(1+\log L_h)$ by comparison with $z\mapsto\alpha_+z+C_0$. Hence, writing
$N_{q,L_h}(Q)=\inf\{n:q^{\circ n}(Q)\le L_h\}$,
\begin{equation}
 h\{N_q(Q_h)-N_{q,L_h}(Q_h)\}\to0.
 \tag{4.2}\label{eq:cutoff}
\end{equation}

Put $c_h=\frac12\log(1/h)$ and
$\nu_h=\inf\{k\ge0:Q_k\le L_h\}$.  While the physical logarithmic depth $Q_k=-\log|X^h_{kh}|$ is above $L_h$,
\begin{equation}
 Q_{k+1}=q(Q_k)+c_h+Y_{k+1}^{h},\qquad
 Y_{k+1}^{h}=-\log|Z_{k+1}+\delta_{h,k}|.
 \tag{4.3}\label{eq:depthrec}
\end{equation}
The affine bounds and stopping at the deep-region crossing keep the shifts $\delta_{h,k}$ uniformly bounded.  Lemma~\ref{lem:logshock} therefore gives, uniformly for those bounded shifts,
\[
 \Pp(Y^h>y)\le C_1e^{-y},\qquad
 \Pp(Y^h<-y)\le C_1e^{-c_1e^{2y}}.
\]
Thus, for each fixed $K$ and any fixed $A>2$,
\begin{equation}
 \Pp\left(\max_{1\le k\le\lfloor K/h\rfloor\wedge\nu_h}
 |Y_k^h|>A\log(1/h)\right)
 \longrightarrow0.
 \tag{4.4}\label{eq:maxshock}
\end{equation}

Let $D_{k+1}=q(D_k)$, $D_0=Q_h$.  Stop the comparison at
$\tau_h^*=\nu_h\wedge N_{q,L_h}(Q_h)\wedge\lfloor K/h\rfloor$.
From \eqref{eq:coarsecontract}, on the high-probability event in \eqref{eq:maxshock}, while $k<\tau_h^*$,
\[
 |Q_{k+1}-D_{k+1}|\le
 \alpha_+|Q_k-D_k|+C+c_h+\max_{j\le K/h}|Y_j^h|.
\]
Consequently, on the event in \eqref{eq:maxshock},
\begin{equation}
 \max_{k\le \tau_h^*}|Q_k-D_k|\le b_h,
 \qquad b_h:=\frac{C+c_h+A\log(1/h)}{1-\alpha_+}
 =O(\log(1/h))=o(L_h).
 \tag{4.5}\label{eq:shadow}
\end{equation}

The $L_h$-threshold is crossed within one further step. Indeed, for $h$ small,
\begin{equation}
 \alpha_+(L_h+b_h)+C_0+c_h+s_h<L_h-b_h,
 \tag{4.6}\label{eq:cross}
\end{equation}
where $s_h=A\log(1/h)$ bounds the shocks on the good event.  If either the stochastic or deterministic orbit first reaches $(-\infty,L_h]$, the other lies below $L_h+b_h$ and \eqref{eq:cross} forces it below $L_h$ at the next step. Therefore,
\begin{equation}
 |\nu_h-N_{q,L_h}(Q_h)|\le1
 \tag{4.7}\label{eq:clockmatch}
\end{equation}
with probability tending to one whenever the relevant crossing occurs before $K/h$.

We now treat the three clock regimes separately.

\emph{Finite positive clock.}  Suppose \(0<d<\infty\), and choose
\(K>d+1\).  Equation \eqref{eq:cutoff} gives
\(hN_{q,L_h}(Q_h)\to d\), so the deterministic crossing occurs before
\(K/h\) for all small \(h\).  The crossing argument then gives \(\nu_h<K/h\) with probability
tending to one.  Hence \eqref{eq:clockmatch} yields
\[
 h\nu_h\to d\quad\hbox{in probability}.
 \tag{4.8}\label{eq:finiteclock}
\]
At $\nu_h$, the good-shock event and the lower affine bound imply that there is a fixed $c_*>0$ such that
\[ c_*L_h\le Q_{\nu_h}\le L_h \]
with probability tending to one. Indeed, $Q_{\nu_h-1}>L_h$ and the lower affine bound in \eqref{eq:affineenvelope} gives
$Q_{\nu_h}\ge\alpha_-Q_{\nu_h-1}-B_-+c_h-s_h\ge\alpha_-L_h-o(L_h)$.  To control the path before the mesoscopic crossing, note that
$Q_k>L_h$ for $k<\nu_h$, so
$|X^h_{kh}|\le e^{-L_h}$.  The lower affine bound in
\eqref{eq:affineenvelope} gives
$\sqrt h\,\sigma(X^h_{kh})\le C\sqrt h\,e^{-\alpha_-L_h}$.
On $\{\nu_h\le K/h\}$, Brownian maximal inequalities yield
\[
\sup_{0\le t\le h\nu_h}|X_t^h|
\le e^{-L_h}+C\sqrt h\,e^{-\alpha_-L_h}
\max_{k<K/h}\sup_{0\le s\le h}
\frac{|W_{kh+s}-W_{kh}|}{\sqrt h}
=o_{\Pp}(1),
\]
since the maximum on the right is
$O_{\Pp}(\sqrt{\log(1/h)})$.
Lemma~\ref{lem:residual} now shows that the remaining
physical time to a fixed rescaled collar is $o_{\Pp}(1)$,
and that the path over this residual interval also vanishes
uniformly in probability.  The rescaled exit position is tight by Lemma~\ref{lem:lyapdescent}.
Let $\tau_h^{\rm out}$ be the first grid-valued exit time from
$(-e^{-B}r_h,e^{-B}r_h)$, with $B$ as in Lemma~\ref{lem:lyapdescent}.
The preceding estimates give
\[
 \tau_h^{\rm out}\xrightarrow{\Pp}d,\qquad
 \sup_{t\le \tau_h^{\rm out}}|X_t^h|\xrightarrow{\Pp}0,
\]
and the rescaled exit states form a tight family. On each fixed rescaled annulus, Proposition~\ref{prop:restart} gives uniform convergence; tightness then permits the annulus to grow after passing to the limit.  For a bounded
continuous functional $G$ of $(\tau_h^{\rm out},X^h_{\cdot\wedge\tau_h^{\rm out}})$
and a bounded Lipschitz functional $F$ of the future path, the strong
Markov property gives
\[
 \E\left[G\,F(X^h_{\tau_h^{\rm out}+\cdot})\right]
 =\E\left[G\,\Psi_h(X^h_{\tau_h^{\rm out}})\right],
 \qquad
 \Psi_h(x)=\E_xF(X^h).
\]
By the uniformity statement in Proposition~\ref{prop:restart} and the tightness of the rescaled
exit state,
\[
 \Psi_h(X^h_{\tau_h^{\rm out}})\xrightarrow{\Pp}\E F(X^{\ns}).
\]
Together with $(\tau_h^{\rm out},X^h_{\cdot\wedge\tau_h^{\rm out}})\to(d,0)$ in
probability, this proves the joint convergence
\[
 (\tau_h^{\rm out},X^h_{\cdot\wedge\tau_h^{\rm out}},X^h_{\tau_h^{\rm out}+\cdot})
 \Rightarrow(d,0,X^{\ns}),
\]
with asymptotic independence of the future component.  The
concatenation map is continuous at this limiting triple because both
pieces are continuous and the pre-restart path is zero.  Hence
$X^h\Rightarrow X^{(d)}$.

\emph{Zero clock.}  If \(d=0\), the affine two-index comparison
\eqref{eq:clockbounds} gives \(N_q(Q)\asymp 1+\log Q\) for large
\(Q\).  Hence \(hN_q(Q_h)\to0\) implies
\(hH(x_h)=h\log(1+Q_h)\to0\), and
Proposition~\ref{prop:restart} applies directly to the original
initialization.  Therefore
\(X^h\Rightarrow X^{\ns}=X^{(0)}\).

\emph{Divergent clock.}  If \(hN_q(Q_h)\to\infty\), fix \(K<\infty\).
Then \eqref{eq:cutoff} gives \(N_{q,L_h}(Q_h)>K/h+1\) eventually.
Were \(\nu_h\le K/h\) on the good shock event, the reverse
crossing implication would give
\(N_{q,L_h}(Q_h)\le K/h+1\), a contradiction.  Thus
\[
 \Pp(\nu_h>K/h)\to1.
 \tag{4.9}\label{eq:infiniteclock}
\]
On \([0,K]\), \eqref{eq:shadow} keeps the grid-point stochastic depth above
\(L_h-O_{\Pp}(\log(1/h))\).  Hence
\[
 \max_{k\le K/h}|X^h_{kh}|\xrightarrow{\Pp}0.
\]
It remains only to control the continuous interpolation.  On the same good-shock event, the lower affine bound for $q$ and the grid-depth lower bound imply uniformly for \(k\le K/h\)
\[
 \sqrt h\,\sigma(X^h_{kh})
 \le \exp\{-cL_h+O(\log(1/h))\}
\]
for some fixed \(c>0\).  Moreover
\(\max_{k\le K/h}\sup_{0\le s\le h}|W_{kh+s}-W_{kh}|/\sqrt h
=O_{\Pp}(\sqrt{\log(1/h)})\) by the Brownian maximal tail and a union bound.  Therefore
\[
 \max_{k\le K/h}\sup_{0\le s\le h}
 \bigl|X^h_{kh+s}-X^h_{kh}\bigr|
 \xrightarrow{\Pp}0,
\]
and consequently
\[
 \sup_{t\le K}|X_t^h|\xrightarrow{\Pp}0.
 \tag{4.10}\label{eq:zeropath}
\]
As \(K\) is arbitrary, the compact-time limit is the zero path.
\end{proof}

\begin{corollary}[Complete cluster-set selection]\label{cor:clusters}
For any deterministic family $x_h\to0$, put
$t_h=hN_q(\log(1/|x_h|))\in[0,\infty]$, with the convention
$t_h=\infty$ when $x_h=0$.  The laws of $X^h$ are relatively compact
in $\mathcal P(C_{\rm loc}([0,\infty)))$, and
\[
 \operatorname{Clust}_{h\downarrow0}\mathcal L(X^h)
 =\{\mathcal L(X^{(d)}):
 d\in\operatorname{Clust}_{h\downarrow0}t_h\},
\]
where the clock cluster set is taken in the compact space
$[0,\infty]$ and $X^{(\infty)}$ denotes the zero path.
\end{corollary}
\begin{proof}
Every sequence of mesh sizes has a subsequence on which $t_h$
converges in $[0,\infty]$.  Theorem~\ref{thm:master} gives convergence
of the associated Euler laws to the corresponding delayed-solution
law.  Conversely, every clock cluster point is attained along a
subsequence, and the same theorem identifies the process-law limit.
\end{proof}

\section{Consequences}
\begin{corollary}[Subcritical equivalence and delay window]\label{cor:window}
Under Assumption~\ref{ass:CI},
\[
 \frac{\log Q}{|\log\alpha_-|}+O(1)
 \le N_q(Q)\le
 \frac{\log Q}{|\log\alpha_+|}+O(1).
 \tag{5.1}\label{eq:windowclock}
\]
Hence $hN_q(Q_h)\to0$ if and only if
$h\log(1+Q_h)\to0$.  If $h\log Q_h\to\theta\in(0,\infty)$, every cluster point $d$ of $hN_q(Q_h)$ belongs to
\[
 \left[\frac{\theta}{|\log\alpha_-|},
       \frac{\theta}{|\log\alpha_+|}\right].
 \tag{5.2}\label{eq:delaywindow}
\]
Every corresponding Euler subsequential limit is $X^{(d)}$. If instead $h\log Q_h\to\infty$, then $hN_q(Q_h)\to\infty$ and the Euler scheme converges to the zero path on compact time intervals.
\end{corollary}
\begin{proof}
Compare the iterates with affine recursions of slopes $\alpha_-$ and $\alpha_+$; additive constants shift the iterates by bounded affine fixed points and affect the hitting count only by $O(1)$ after taking logarithms.  Theorem~\ref{thm:master} gives the final assertion.
\end{proof}

\begin{corollary}[Common-index generalized power law]\label{cor:GP}
Suppose that for some $0<\alpha<1/2$ and every $\eta>0$ there are $C_\eta,\delta_\eta>0$ such that
\[
 C_\eta^{-1}t^{\alpha+\eta}\le
 \frac{\sigma(tx)}{\sigma(x)}\le
 C_\eta t^{\alpha-\eta},\qquad0<t\le1,\quad0<x<\delta_\eta.
 \tag{GP}
\]
Then
\[
 N_q(Q)\sim\frac{\log Q}{|\log\alpha|}.
\]
Consequently, if
\[
 h\log\!\left(1+\log^+\frac1{|x_h|}\right)\to\theta\in[0,\infty],
\]
then
\[
 X^h\Rightarrow X^{(\theta/|\log\alpha|)}.
\]
For every $d>0$, the initialization
\[
 x_h=\exp\{-\exp(d|\log\alpha|/h)\}
\]
realizes the delay $d$ for every coefficient satisfying (GP).
\end{corollary}
\begin{proof}
For each $\eta>0$, (GP) gives \eqref{eq:CI} with indices $\alpha-\eta$ and $\alpha+\eta$ and an additive constant $\log C_\eta$. Apply \eqref{eq:windowclock} and let $\eta\downarrow0$. The displayed initialization has $h\log\log(1/|x_h|)=d|\log\alpha|$.
\end{proof}

Regularly varying coefficients $\sigma(x)=|x|^\alpha\ell(|x|)$ are included, by Potter's bounds \cite{BGT1987}.  The next result shows that a single effective clock rate may survive even when no asymptotic power exponent exists.

\begin{proposition}[Periodic logarithmic profiles]\label{prop:periodic}
Let $m\in C^1(\mathbb R)$ be positive and $P$-periodic, and suppose that
\[
 0<\alpha_-\le m(s)+m'(s)\le\alpha_+<\frac12,
 \qquad s\in\mathbb R.
 \tag{5.3}\label{eq:periodicass}
\]
Suppose that
\begin{equation}
 q(Q)=Qm(\log Q)
 \tag{5.4}\label{eq:periodicq}
\end{equation}
for all sufficiently large $Q$.
Then \eqref{eq:CI} holds with additive constant zero.  Moreover, with
\[
 F(s)=s+\log m(s),
 \tag{5.5}\label{eq:Fperiodic}
\]
there exists a number $\rho<0$ such that
\[
 \sup_{s\in\mathbb R}\left|F^{\circ n}(s)-s-n\rho\right|\le P,
 \qquad n\ge1,
 \tag{5.6}\label{eq:boundeddeviation}
\]
and the fixed-reference iteration clock satisfies
\[
 N_q(Q)=\frac{\log Q}{|\rho|}+O(1),\qquad Q\to\infty.
 \tag{5.7}\label{eq:periodicclock}
\]
In addition,
\[
 \log\alpha_-\le\rho\le\log\alpha_+<0.
 \tag{5.8}\label{eq:rhobounds}
\]
\end{proposition}
\begin{proof}
Since $q'(Q)=m(\log Q)+m'(\log Q)$, \eqref{eq:periodicass} gives \eqref{eq:CI} with $C=0$.  At a minimum and a maximum of the periodic function $m$ one has $m'=0$; hence
$\alpha_-\le m(s)\le\alpha_+$ for every $s$.  Thus
\[
 \log\alpha_-\le F(s)-s\le\log\alpha_+<0.
 \tag{5.9}\label{eq:Fstep}
\]
Also
\[
 F(s+P)=F(s)+P,
 \qquad
 F'(s)=\frac{m(s)+m'(s)}{m(s)}>0.
\]
Hence $F$ is an increasing lift of an orientation-preserving circle homeomorphism of circumference $P$.

Set
\[
 A_n=\sup_s\{F^{\circ n}(s)-s\},\qquad
 B_n=\inf_s\{F^{\circ n}(s)-s\}.
\]
Monotonicity and $P$-equivariance imply $A_n-B_n\le P$.  Moreover,
$A_{n+k}\le A_n+A_k$ and $B_{n+k}\ge B_n+B_k$.  The standard translation-number argument for lifts of orientation-preserving circle homeomorphisms (see, e.g., \cite[Ch.~11]{KatokHasselblatt1995}) therefore yields a common $\rho$ with
$B_n\le n\rho\le A_n$; consequently \eqref{eq:boundeddeviation} holds.  Summing \eqref{eq:Fstep} along an orbit gives \eqref{eq:rhobounds}.

Choose a fixed $Q_{\rm per}$ so large that \eqref{eq:periodicq} holds for every $Q\ge Q_{\rm per}$, and put $S_{\rm per}=\log Q_{\rm per}$.  Let
\[
 N_{\rm per}(Q):=\inf\{n\ge0:q^{\circ n}(Q)\le Q_{\rm per}\}.
\]
For every $n<N_{\rm per}(Q)$, with $S=\log Q$,
$q^{\circ n}(Q)=\exp(F^{\circ n}(S))$.  At the first crossing of $S_{\rm per}$, the two defining inequalities, together with \eqref{eq:boundeddeviation} and $\rho<0$, give
\[
 N_{\rm per}(Q)=\frac{S-S_{\rm per}}{|\rho|}+O(1)
 =\frac{\log Q}{|\rho|}+O(1).
\]
Both $Q_{\rm per}$ and $Q_{\rm ref}$ are fixed admissible reference levels after enlarging $Q_{\rm per}$ if necessary.  Lemma~\ref{lem:detclock} therefore gives
$N_q(Q)=N_{\rm per}(Q)+O(1)$, which proves \eqref{eq:periodicclock}.
\end{proof}

\begin{corollary}[Effective delay for periodic logarithmic modulation]\label{cor:periodic}
Under the assumptions of Proposition~\ref{prop:periodic}, let $x_h\to0$ be deterministic and put $Q_h=\log(1/|x_h|)$. If
\[
 h\log Q_h\longrightarrow\theta\in[0,\infty),
\]
then
\[
 X^h\Rightarrow X^{(\theta/|\rho|)}.
 \tag{5.10}\label{eq:periodicdelay}
\]
In particular, for $m(s)=a+b\sin s$, condition \eqref{eq:periodicass} is satisfied whenever
\[
 0<a-\sqrt2|b|,
 \qquad
 a+\sqrt2|b|<\frac12.
 \tag{5.11}\label{eq:sincond}
\]
Near zero the corresponding coefficient is
\[
 \sigma(x)=|x|^{\,a+b\sin(\log\log(1/|x|))}.
 \tag{5.12}\label{eq:oscillatorysigma}
\]
Its local power exponent has no limit when $b\ne0$, while the Euler selection clock has the single effective rate $1/|\rho|$.
\end{corollary}
\begin{proof}
Equation \eqref{eq:periodicclock} gives
$hN_q(Q_h)=h\log Q_h/|\rho|+O(h)$, and Theorem~\ref{thm:master} applies.  For $m(s)=a+b\sin s$,
$m(s)+m'(s)=a+b(\sin s+\cos s)$, so \eqref{eq:sincond} implies \eqref{eq:periodicass}.  Formula \eqref{eq:oscillatorysigma} follows from
$\sigma(e^{-Q})=e^{-q(Q)}$.
\end{proof}

The number $\rho$ is the translation number of the real lift $F$, not a rotation number modulo $P$. For constant $m\equiv\alpha$, it equals $\log\alpha$ and Corollary~\ref{cor:periodic} gives the familiar power-law delay $\theta/|\log\alpha|$.

\section{Full sharpness of the two-index delay window}
Fix $0<\alpha_-<\alpha_+<1/2$ and $\theta>0$, and write
\[
 I_\theta=\left[\frac{\theta}{|\log\alpha_-|},
                    \frac{\theta}{|\log\alpha_+|}\right].
\]
We construct one coefficient and one initialization family whose
complete clock cluster set is exactly $I_\theta$.  Here $\operatorname{Clust}_{h\downarrow0}$ denotes the set of cluster points along sequences $h\downarrow0$.

Let $D=\{d_1,d_2,\ldots\}$ be a countable dense subset of $I_\theta$,
and choose a sequence $\widehat d_m$ in which every $d_j$ occurs
infinitely often.  Put
\[
 \beta_m=\exp(-\theta/\widehat d_m)\in[\alpha_-,\alpha_+].
 \tag{6.1}\label{eq:betam}
\]
Construct a continuous slope field
$a:[Q_\star,\infty)\to[\alpha_-,\alpha_+]$ recursively and set
\[
 q(Q)=q_0+\int_{Q_\star}^Q a(s)\,ds.
 \tag{6.2}\label{eq:qintegral}
\]
Then globally
\[
 \alpha_-(Q_2-Q_1)\le q(Q_2)-q(Q_1)
 \le\alpha_+(Q_2-Q_1),
 \tag{6.3}\label{eq:sharpCI}
\]
so (CI) holds with additive constant zero.

At stage $m$, after all earlier stages have been completed, we move to an
arbitrarily large new depth and connect continuously to the slope
$\beta_m$.  Continue with $a\equiv\beta_m$ on a pure affine block.
Choose $R_m$ arbitrarily far inside this block and define
\[
 h_m=\exp\{-2(R_m-q(R_m))\},\qquad
 Q_{0,m}=e^{\theta/h_m}.
 \tag{6.4}\label{eq:hm}
\]
By taking $R_m$ farther to the right if necessary, arrange
$h_m<\min\{1/m,h_{m-1}/2\}$ (with the second constraint omitted for
$m=1$), and then extend the same pure block through $Q_{0,m}+1$.
All later stages are constructed strictly to the right of this point.
All iterates of $q$ from $Q_{0,m}$ down to $L_m$ therefore lie in the stage-$m$ affine block.

On the stage-$m$ block,
$q(Q)=\beta_mQ+\gamma_m$.  With
$Q_{*,m}=\gamma_m/(1-\beta_m)$, take the lower block level
$L_m=R_m$.  By choosing $R_m\uparrow\infty$ recursively, we have
$[L_m,Q_{0,m}]$ inside the same pure block and, because
$h_m=\exp\{-2(R_m-q(R_m))\}$ with
$R_m-q(R_m)\ge(1-\alpha_+)R_m-O(1)$,
also $h_m\log L_m\to0$.  Since $R_m-q(R_m)=\tfrac12\log(1/h_m)$ and
$q(Q)=\beta_mQ+\gamma_m$ on the block,
\[
 L_m-Q_{*,m}=R_m-Q_{*,m}
 =\frac{\log(1/h_m)}{2(1-\beta_m)}>0,
\]
so $h_m\log(L_m-Q_{*,m})\to0$.  The exact affine iteration formula gives
\[
 q^{\circ n}(Q)-Q_{*,m}=\beta_m^n(Q-Q_{*,m}),
\]
and therefore the number $N_m$ of iterates from $Q_{0,m}$ to $L_m$
satisfies
\[
 h_mN_m\longrightarrow
 \frac{\theta}{|\log\beta_m|}=\widehat d_m.
 \tag{6.5}\label{eq:localsharpclock}
\]
The residual number of iterates from $L_m$ to the fixed reference
level is $O(\log L_m)$ by Lemma~\ref{lem:detclock}, hence negligible
after multiplication by $h_m$.  Thus
\[
 h_mN_q(Q_{0,m})-\widehat d_m\longrightarrow0.
 \tag{6.6}\label{eq:sharpclock}
\]

Define near zero
\[
 \sigma(x)=\exp\{-q(\log(1/|x|))\},\qquad x\ne0,
 \qquad \sigma(0)=0,
\]
and extend it continuously and positively away from zero with at most
linear growth.  This gives one admissible coefficient.  Because every
$d_j$ occurs infinitely often and the construction already enforces
$h_m\downarrow0$, each $d_j$ is a cluster point of
$hN_q(e^{\theta/h})$.  Corollary~\ref{cor:window} gives the reverse
inclusion in $I_\theta$, and the cluster set is closed.  Hence
\[
 \operatorname{Clust}_{h\downarrow0}hN_q(e^{\theta/h})=I_\theta.
 \tag{6.7}\label{eq:fullclockcluster}
\]

\begin{theorem}[Full sharpness of the delay window]\label{thm:sharp}
There exists a coefficient $\sigma$ satisfying \eqref{eq:CI} with $C=0$ and $\alpha_-<\alpha_+$ such that,
for the deterministic initialization family
\[
 x_h=\exp\{-\exp(\theta/h)\},
\]
we have
\[
 \operatorname{Clust}_{h\downarrow0}hN_q(\log(1/|x_h|))=I_\theta.
\]
Moreover the cluster set of the Euler laws in
$\mathcal P(C_{\rm loc}([0,\infty)))$ is exactly
\[
 \{\mathcal L(X^{(d)}):d\in I_\theta\}.
\]
Thus the delay window \eqref{eq:delaywindow} is optimal for a single
coefficient and initialization family.
\end{theorem}
\begin{proof}
The clock statement follows from \eqref{eq:fullclockcluster}.
The claim for the process laws follows from
Corollary~\ref{cor:clusters}.
\end{proof}

\section{Discussion}
For the homogeneous power coefficient studied in the companion paper, the clock reduces to the explicit double-logarithmic selector.  The general theorem instead describes selection through the coefficient-dependent passage-time clock.
In the common-index regime the iteration count has the closed form
$N_q(Q)\sim\log Q/|\log\alpha|$.  Periodic logarithmic profiles show that a single effective rate may persist even without an asymptotic power exponent: it is then the reciprocal absolute translation speed of the induced logarithmic dynamics.  With general separated indices the clock may retain still finer coefficient-specific information, constrained by the delay window \eqref{eq:delaywindow}.  The sharpness construction realizes the entire two-index window with one coefficient and one initialization family.  The upper-index condition $\alpha_+<1/2$ ensures local integrability of $\sigma^{-2}$ and the vanishing occupation estimates used to identify the fundamental solution.  For the power coefficient $\sigma(x)=|x|^\alpha$, the case $\alpha\ge1/2$ has the absorbing solution as its unique weak solution from zero; hence the threshold is structural for that family. For general coefficients the upper-index condition is sufficient, not necessary. In particular, the class might be enlarged under assumptions formulated directly in terms of
$J(r)=\int_0^r\sigma(a)^{-2}\,da$.  For example, $\sigma(x)=\sqrt{|x|}\log(1/|x|)$ still has locally integrable $\sigma^{-2}$ despite its nominal power exponent $1/2$.  Establishing an iteration-clock selection theorem in such a borderline class would require new occupation estimates. The local symmetry assumption permits a single depth map $q$ for both signs. For asymmetric coefficients, the signs of successive Euler updates can select different logarithmic profiles, and the present deterministic clock no longer directly describes the recursion.

The finite-delay regime is asymptotic: it requires logarithmic depths of order $\exp(c/h)$, so the corresponding physical initializations eventually underflow in any fixed finite-precision representation.  Numerical exploration of the critical regime is therefore more naturally carried out through the logarithmic-depth recursion than through the physical initialization itself.  Further directions include second-order fluctuations of the clock, asymmetric coefficients, and comparison of numerical schemes through their induced logarithmic iteration maps.

\section*{Disclosure}
The author used OpenAI ChatGPT to assist with language editing, manuscript organization, literature searches, and discussion and critical examination of mathematical arguments and proof strategies. The author reviewed the manuscript, made the decisions regarding its content, and takes full responsibility for its originality, correctness, and integrity.

\par\medskip
\noindent The views expressed in this article are solely those of the author and do not necessarily reflect those of Sumitomo Mitsui Banking Corporation.

\end{document}